\documentclass{article}    
\usepackage{latexsym}
\usepackage{amsbsy}
\usepackage{amssymb}
\usepackage[latin1]{inputenc} 

\title{Diffeomorphism invariant Colombeau algebras. Part I: Local theory}
\author{Roland Steinbauer\footnote{Electronic mail: roland.steinbauer@univie.ac.at;
supported by research grant P12023-MAT of the Austrian Science Foundation} \\
  \\{\small Department of Mathematics, University of Vienna}\\
  {\small Strudlhofg. 4, A-1090 Wien, Austria}}
\date{}
\begin{document}
\newcommand{\al}{\alpha}\newcommand{\be}{\beta} 
\newcommand{\bet}{\beta}\newcommand{\ga}{\gamma}
\newcommand{\Om}{\Omega}\newcommand{\Ga}{\Gamma}\newcommand{\om}{\omega}
\newcommand{\si}{\sigma}\newcommand{\la}{\lambda}
\newcommand{\eps}{\varepsilon}\newcommand{\de}{\delta}
\newcommand{\vphi}{\varphi}\newcommand{\dl}{{\displaystyle \lim_{\eta>0}}\,}
\newcommand{\intl}{\int\limits}\newcommand{\su}{\sum\limits_{i=1}^2}
\newcommand{\C}{{\mathbb C}}\newcommand{\R}{{\mathbb R}}\newcommand{\N}{{\mathbb N}}
\newcommand{\K}{{\mathbb K}}
\newcommand{\D}{{\cal D}}\newcommand{\Vol}{\mathrm{Vol\,}}
\newcommand{\tr}{\mathrm{tr}}\newcommand{\spn}{\mathrm{span}}
\newcommand{\Or}{\mbox{Or}}\newcommand{\sign}{\mbox{sign}}
\newcommand{\na}{\nabla}\newcommand{\pa}{\partial}
\newcommand{\ti}{\tilde}\newcommand{\T}{{\cal T}} \newcommand{\G}{{\cal G}}
\newcommand{\DD}{{\cal D}}\newcommand{\SSS}{{\cal S}}
\newcommand{\X}{{\cal X}}\newcommand{\E}{{\cal E}} 
\newcommand{\CC}{{\cal C}}\newcommand{\vo}{\Vol}
\newcommand{\NN}{{\cal N}}
\newcommand{\bt}{\bar t}\newcommand{\bx}{\bar x}
\newcommand{\by}{\bar y} \newcommand{\bz}{\bar z}\newcommand{\br}{\bar r}
\newcommand{\fr}{\frac{1}}\newcommand{\il}{\int\limits}
\newcommand{\nn}{\nonumber}\newcommand{\supp}{\mathrm{supp}}
\newcommand{\vp}{\mbox{vp}\frac{1}{x}}\newcommand{\A}{{\cal A}}
\newcommand{\Ll}{L_{\mbox{\small loc}}}\newcommand{\Hl}{H_{\mbox{\small loc}}}
\newcommand{\Lll}{L_{\mbox{\scriptsize loc}}}
\newcommand{\cl}{\mbox{\rm cl}}
\newcommand{\gK}{{\cal K}}\newcommand{\gR}{{\cal R}}\newcommand{\gC}{{\cal C}}
\newcommand{\ddet}{\mathop{\mathrm{det}}}
\newcommand{\lla}{\langle}\newcommand{\ra}{\rangle}
\newcommand{\beq}{ \begin{equation} }\newcommand{\eeq}{\end{equation} }
\newcommand{\bea}{\begin{eqnarray}}\newcommand{\eea}{\end{eqnarray}}
\newcommand{\beas}{\begin{eqnarray*}}\newcommand{\eeas}{\end{eqnarray*}}
\newcommand{\beqs}{\begin{equation*}}\newcommand{\eeqs}{\end{equation*}}
\newcommand{\lb}{\label}\newcommand{\rf}{\ref}
\newcommand{\GL}{\mbox{GL}}\newcommand{\bfs}{\boldsymbol}
\newcommand{\ben}{\begin{enumerate}}\newcommand{\een}{\end{enumerate}}
\newcommand{\ba}{\begin{array}}\newcommand{\ea}{\end{array}}
\newtheorem{thi}{\hspace*{-1.1mm}}[section]
\newtheorem{thr}{thi}
\newcommand{\bthm}{\begin{thi} {\bf Theorem. }}
\newcommand{\bprop}{\begin{thi} {\bf Proposition. }}
\newcommand{\bcor}{\begin{thi} {\bf Corollary. }}
\newcommand{\blem}{\begin{thi} {\bf Lemma. }}
\newcommand{\bex}{\begin{thi} {\bf Example. }\rm}
\newcommand{\bexs}{\begin{thi} {\bf Examples. }\rm}
\newcommand{\bdef}{\begin{thi} {\bf Definition. }}
\newcommand{\brem}{\begin{thi} {\bf Remark. }\rm}
\newcommand{\bth}{\begin{thi}\rm}
\newcommand{\ethi}{\end{thi}}
\newcommand{\bblock}{\\ \begin{minipage}{15cm} \beas}
\newcommand{\eblock}[1]{\eeas \end{minipage} 
\hfill \begin{minipage}{10mm}\bea \label{#1}\eea\end{minipage}\\}

\newcommand{\ep}{\hspace*{\fill}$\Box$}
\newcommand{\ms}{\medskip\\}
\newcommand{\et}{\end{thi}}
\newcommand{\lgl}{\langle}
\newcommand{\rgl}{\rangle}
\newcommand{\pr}{{\bf Proof. }}
\newcommand{\cc}{{\cal C}}
\newcommand{\comp}{\subset\subset}
\newcommand{\ca}{{\cal A}}
\newcommand{\cb}{{\cal B}}
\newcommand{\cd}{{\cal D}}
\newcommand{\ce}{{\cal E}}
\newcommand{\cg}{{\cal G}}
\newcommand{\ci}{{\cal I}}
\newcommand{\cn}{{\cal N}}
\newcommand{\cs}{{\cal S}}
\newcommand{\rmd}{\mbox{\rm d}}
\newcommand{\io}{\iota}
\newcommand{\bnot}{\begin{thr} {\bf Notation }}
\newcommand{\spp}{\mbox{\rm supp\,}}
\newcommand{\id}{\mathop{\mathrm{id}}}
\newcommand{\pro}{\mathop{\mathrm{pr}}}
\newcommand{\dist}{\mathop{\mathrm{dist}}}
\newcommand{\clb}{\overline{B}_}
\newcommand{\ahat}{\hat{\mathcal{A}}_0(X)} 
\newcommand{\atil}{\tilde{\mathcal{A}}_0(X)} 
\newcommand{\ehat}{\hat{\mathcal{E}}(X)} 
\newcommand{\emhat}{\hat{\mathcal{E}}_M(X)} 
\newcommand{\nhat}{\hat{\mathcal{N}}(X)}
\newcommand{\cee}{{\cal E}}
\newcommand{\todo}[1]{$\clubsuit$\ {\tt #1}\ $\clubsuit$} 
\maketitle
\begin{abstract}     
This contribution is the first in a series of three: it reports on the construction
of (a fine sheaf of) diffeomorphism invariant Colombeau algebras $\G^d(\Om)=
\E_M(\Om)/\NN(\Om)$ on open sets of Eucildean space (\cite{found}), 
which completes earlier approaches (\cite{CM,Jel}). Part II and III will 
show, among others, the way to an intrinsic definition of Colombeau algebras on 
manifolds which, locally, reproduces the algebra(s) $\G^d(\Om)$.
\end{abstract}

\section{Introduction}
   
Since its introduction in \cite{c1} the question of the functor property of
Co\-lom\-beau's construction was at hand as a crucial one:
let $\mu:\ti\Om\to\Om$ be a diffeomorphism between open subsets of $\R^s$, is
it possible to extend the operation $\mu^*:f\mapsto f\circ\mu$ on
smooth distributions on $\Om$ to an operation $\hat\mu$ on
the Colombeau algebra such that $(\mu\circ\nu)\hat{\ }=\hat\nu\circ\hat\mu$ and 
$(\mbox{\rm id})\hat{\ }=\mbox{\rm id}$ are satisfied or---to phrase it differently---is 
it possible to achieve a diffeomorphism invariant construction of {\em full} Colombeau 
algebras? (In this work we shall focus on full algebras distinguished by the existence 
of a canonical embedding of distributions and henceforth omit the term ``full''. 
By contrast the so called {\em special} 
(or simplified) algebras with their elements basically depending on a real 
regularization parameter $\eps\in(0,1]$ do not allow for a canonical embedding of 
distributions. However, by their relative ease of construction and the fact that 
diffeomorphism invariance of the basic definitions is automatically satisfied 
they provide a flexible tool to model singularities in an nonlinear context: global
analysis in this setting has been investigated in \cite{ndg}.)

To begin with we briefly recall Colombeau's construction as given in \cite{c2} 
and shall refer to the corresponding algebra as the ``elementary'' one. 
Let $\Om\subseteq\R^s$ and for $\vphi\in\D(\R^s)$ set $\vphi_\eps(x)=
(1/\eps)^s\vphi(x/\eps)$. Then we define:
\beas 
 {\cal  A}_0(\R^s)  &:=&  \{\vphi  \in  {\cal D}(\R^s) \, : \, \int 
\vphi(x)\,dx = 1\}  \\ 
 {\cal  A}_q(\R^s)  &:=&  \{\vphi  \in  {\cal  A}_0(\R^s) \, : \, \int 
\vphi(x)  x^\alpha  \,dx  =  0\,  ,  \,  1\le  |\alpha| \le q\} 
\,\,\,(q \in \N)\\ 
U(\Om)&:=&\{(\vphi,x)\in{\cal A}_0(\R^s)\times\Om\,:\,\supp(\vphi)\subseteq\Om-x\}\\
{\cal E}^e(\Omega) &:=& \{R: U(\Omega) \to \C \, :\,\forall (\vphi,x)\in U(\Om)
\mbox{ the map }
y\mapsto R(\vphi,y)\mbox{ is}\\&&\mbox{ smooth in a neighborhood of } x\,\} \\ 
{\cal E}^e_M(\Omega) &:=& \{R\in {\cal E}(\Omega) :  
  \forall K\subset\subset \Omega \ \forall \alpha\in \N_0^s \ \exists N\in 
  \N_0 \ \forall \varphi\in {\cal A}_p(\R^s):\\&&\ \sup\limits_{x\in K} | \partial^\alpha 
  R(\varphi_\eps,x) | =O(\eps^{-N}) \mbox{ as }\eps\to 0\}  \\
{\cal N}^e(\Omega) &:=& \{R\in {\cal E}(\Omega) :  
  \forall K\subset\subset \Omega \ \forall \alpha\in \N_0^s \ \forall n\in\N_0 \ 
\exists q \ \forall \varphi\in {\cal A}_q(\R^s):\\ &&\ 
  \sup\limits_{x\in K} | \partial^\alpha 
  R(\varphi_\eps,x) | =O(\eps^n)\mbox{ as }\eps\to 0\}\,.
\eeas 
The {\em ``elementary'' Colombeau algebra} of generalized functions on 
$\Om$ finally is defined as the quotient space 
\[
        {\cal G}^e(\Omega):={\cal E}^e_M(\Omega)\,/\,{\cal N}^e(\Omega)\,.
\]

$\G^e(\_)$ is a fine sheaf of differential algebras. Smooth functions are embedded
into $\G^e$ by the ``constant'' embedding $\sigma$, i.e., $\sigma(f)(\vphi,x)=f(x)$, 
turning $\CC^\infty$ into a faithful subalgebra. Distributions,
on the other hand, are embedded via (anti)convolution with the $\vphi$'s, i.e.
$\iota:\,{\cal D}'(\Om)\to{\cal G}^e(\Om)$ with
\begin{equation}\label{iota}
        \iota(u)(\vphi,x)\,=\,\lgl u,\vphi(.-x)\rgl\,.
\end{equation}
It is one of the fundamental properties of Colombeau algebras that $\iota|_{
\CC^\infty}=\sigma$.

Now, let $\mu:\tilde\Om\to\Om$ denote a diffeomorphism as above.
Accepting $\lgl\mu^*u,\vphi\rgl:=\lgl u,(\vphi\circ\mu^{-1})\cdot|\det D\mu^{-1}|\rgl$
and (\ref{iota}) as the explicit forms of $\mu^*:\cd'(\Om)\to\cd'(\ti\Om)$
resp.\ $\io:\cd'\to\cee$, the basic requirement $\hat\mu\circ\io=\io\circ\mu^*$ 
amounts to
$$
(\hat\mu(\io u))(\ti\vphi,\ti x)=
   (\io u)(\left(\ti\vphi(\mu^{-1}(.+\mu\ti x)-\ti x)\cdot
        |\det D\mu^{-1}(.+\mu\ti x)|\,,\,\mu\ti x\right)
$$
(where $(\ti\vphi,\ti x)\in U(\ti\Om)$), which in turn, enforces
\beq\label{mutrsf}
   (\hat\mu R)(\ti\vphi,\ti x):=
   R(\left(\ti\vphi(\mu^{-1}(.+\mu\ti x)-\ti x)\cdot
        |\det D\mu^{-1}(.+\mu\ti x)|\,,\,\mu\ti x\right)
\eeq
as the definition of the action of the diffeomorphism on a representative $R$ 
of a generalized function in a Colombeau algebra on $\Om$.


\section{Towards diffeomorphism invariance}

Colombeau and Meril in their paper \cite{CM} (using earlier ideas 
of \cite{c1}) made the first decisive steps to incorporate formula (\ref{mutrsf}) 
into the construction of a Colombeau algebra which they claimed to be diffeomorphism 
invariant. Before discussing their work in some more detail let us introduce some 
terminology (cf.\ \cite{found}, Section 9) which eases the understanding of the 
definitions to be given below.

Every  Colombeau algebra is constructed as a quotient space
of moderate modulo negligible sequences (nets) of (smooth) functions
$R$ belonging to some {\em basic space} usually denoted by $\E$ (plus
some superscript to distinguish the algebras to be constructed).  
The respective properties of moderateness and negligibility are then defined by
inserting {\em scaled test objects} (e.g.\ $\vphi_\eps$ with
$\vphi\in\A_0$ resp.\ $\A_q$ as above) into $R$ and analyzing the asymptotic 
behavior of the latter on these ``paths'' as the scaling parameter $\eps$ tends to 
zero (and consequently $\vphi_\eps\to\de$ weakly):
we shall refer to the respective processes as {\em testing for moderateness resp.\
negligibility}. In this terminology diffeomorphism invariance of a
Colombeau algebra is ensured by diffeomorphism invariance of the
respective tests, of course including diffeomorphism invariance of the
respective class of (scaled) test objects.
As opposed to this testing procedure the elements of the algebra themselves 
{\em do not depend in any way from $\eps$}. We regard this distinction as 
fundamental clarifying several misinterpretations in the literature and call it 
the policy of {\em ``separating definitions from testing.''} 

In \cite{found}, Section 3 there was given a blueprint collecting all the definitions
and theorems necessary for the construction of a  Colombeau algebra. In the 
following we shall use this collection as a guiding line in discussing the various 
variants of the algebra proposed in the literature beginning with the one of Colombeau 
and Meril \cite{CM}.

There are basically three modifications introduced by the authors of \cite{CM}
distinguishing their construction---which we call $\G^1$---from $\G^e$, namely:
\begin{itemize}
\item[(i)] Smooth dependence of $R$ on $\varphi$ in place of arbitrary dependence.
\item[(ii)] Dependence of test objects on $\eps$, i.e.,  bounded paths 
$\eps\mapsto\phi(\eps)\in\D(\Om)$.
\item[(iii)] Asymptotically vanishing moments (see below) of test objects as compared 
to the stronger condition $\phi(\eps)\in \A_q(\R^s)$ for all $\eps$ (which is the
naive analog of $\vphi\in\A_q(\R^s)$ in the case of $\G^e$).
\end{itemize}
Condition (i) is necessary to guarantee smoothness of $\hat\mu R$ with respect 
to $\tilde x$ (cf.\ transformation (\ref{mutrsf})). However, the technical prize 
to pay here is the use of calculus in infinite dimensional spaces: 
Colombeau and Meril in particular used the concept of Silva-differentiability 
\cite{c0}. However, instead of giving the proofs they rather ``invited the 
reader to admit''(\cite{CM}, p. 362) the respective smoothness properties.

Change (ii) together with (iii) obviously was introduced to obtain a diffeomorphism 
invariant analog of the vanishing moment conditions defined above. 
More precisely, define  
$\widetilde{{\cal A}}_q(\R^s)$ to be the set of all smooth, bounded paths
$\eps \mapsto \phi(\eps)$ satisfying
\bea
        \int \phi(\eps)(\xi)\,d\xi&=& 1 \quad \forall \eps \in (0,1]\mbox{ and}\nn\\
        \label{ColM}
        \int x^\alpha \phi(\eps)(\xi)\, d\xi &=& O(\eps^q) 
        \quad \forall \alpha\in \N_0^s \mbox{ with } 1<|\alpha| \le q\,.
\eea
It may now be shown (\cite{CM}, \S 3) that these moment conditions 
indeed are invariant under the action of a diffeomorphism.

Colombeau and Meril chose their basic space to be $\E^1(\Om):=
\{R:\,\widetilde{{\cal A}}_0(\R^s)\times\Om\to\C^{(0,1]}\}$. Note that this definition
is not in accordance with
the policy of ``separating definitions from testing'' as propagated above.
Moreover, their definition of the objects constituting the Colombeau algebra was not 
unambiguous. However, following the interpretations of \cite{Jel} and \cite{found}, 
the testing process in \cite{CM} is defined by inserting test objects of the form
$S_\eps\phi(\eps):=(1/\eps)^s\vphi(\eps)(./\eps)$ into the first slot of $R$.
More precisely, 
\beas
        {\cal E}^1_M(\Omega) 
        &:=& \{R\in {\cal E}^1(\Omega) :\ \forall 
                K\subset\subset\Omega\ \forall\alpha\in \N_0^s\ \exists N\in \N 
        \forall \phi\in \widetilde{{\cal A}}_p: \\ 
        &&\ \sup\limits_{x\in K}|\partial^\alpha R(S_\eps\phi(\eps),x)|= O(\eps^{-N})
        \mbox{ as } (\eps \to 0)\}\\
        {\cal N}^1(\Omega)
        &:=& \{R\in {\cal E}^1(\Omega) :
        \ \forall 
        K \subset\subset \Omega \ \forall \alpha\in \N_0^s\ \forall n\in \N 
        \ \exists q \ \forall \phi\in \widetilde{{\cal A}}_q: \\
        &&\ \sup\limits_{x\in K}|\partial^\alpha R(S_\eps\phi(\eps),x)|= O(\eps^n)\}
\eeas
Finally the {\em Colombeau-Meril algebra } on $\Om$ is defined as the quotient space
\[
        {\cal G}^1(\Omega) := {\cal E}^1_M(\Omega)\,/\,
        {\cal N}^1(\Omega)\,. 
\]
Using these definitions, all the main properties of ${\cal G}^e(\Omega)$
carry over to ${\cal G}^1(\Omega)$, with almost identical proofs. Indeed, 
boundedness of the paths $\phi(\eps)$ in ${\cal D}(\R^s)$ assures similar 
estimates as in the case of single functions $\varphi$.

Unfortunately, in addition to the ambiguities mentioned above the class of test 
objects as defined by Colombeau and Meril still is not preserved under the action of a 
diffeomorphism. Nevertheless, despite these defects (which, apparently, went unnoticed 
by nearly all workers in the field) their construction was quoted and used many times
(among others \cite{kunzDISS}, \cite{invc}, \cite{NPS}, \cite{geo2}).
It was only in 1998 that J. Jel\'\i nek in \cite{Jel} pointed out the
error in \cite{CM} by giving a (rather simple) counter example which we shall discuss
in a moment.  In the same paper, he presented another version of the
theory which avoided (some of) the shortcomings of \cite{CM} and has
to be considered as the second decisive step towards a diffeomorphism invariant 
version of a Colombeau algebra.

Taking a closer look on the nature of test objects as used by Colombeau and Meril,
from (\ref{mutrsf}) we see that the action of a diffeomorphism $\mu$ introduces
an additional $x$-dependence in the first slot of $R$. This in turn may be 
exploited by giving an example of a function in $\E^1_M$ which is constant in 
$x$ (hence the estimates of the derivatives follow trivially) but whose $\mu$-transform 
depending on $x$ fails to be moderate. More precisely,
set $R(\phi,x):=\exp(i\exp(\int|\phi(\xi)|^2\,d\xi))$.
Then according to \ref{mutrsf} we have 
\beas
        \hat\mu R(S_\eps\tilde\phi,\tilde x)
        &=&R(\bar\mu(S_\eps\tilde\phi,\tilde x))\\
        &=&\exp(i\exp(\int|\tilde\phi(\frac{\mu^{-1}(\eps\xi+\mu\tilde x)-\tilde
                x}{\eps})\det D\mu^{-1}(\eps\xi+\mu\tilde x)|^2\,d\xi))\,.
\eeas

We next discuss in some detail the algebra proposed by J. Jel\'\i nek in \cite{Jel} 
which we shall call $\G^d$ \footnote{``$d$'' obviously
stands for diffeomorphism invariant. In fact Jel\'\i nek's construction comes
very close to diffeomorphism invariance although the last gaps were only closed
in \cite{found}, however, without the necessity to change the main definitions.}. 
Analogously to the previous construction we
start by listing the main features distinguishing $\G^d$ from $\G^1$.
\begin{itemize}
\item[(i)] (Smooth) dependence of test objects also on $x\in\Om$.
\item[(ii)] In testing for moderateness test objects may take arbitrary values
        in $\A_0(\R^s)$, independently of any moment condition.
\end{itemize}
While in the light of Jel\'\i nek's counterexample (i) is compelling
there seems to be no apparent necessity for (ii). Apparently (ii) widens
the range of test objects which in turn reduces $\E^d_M$ resp.\ $\NN^d$
in size. Yet is has to be admitted that by this reduction no generalized function
of interest, neither for the further development of the theory nor in applications
is lost. For the construction of a diffeomorphism invariant Colombeau algebra omitting
(ii) see Part II resp.\ \cite{found}. Here, however, we focus on $\G^d$ which we regard 
as the standard diffeomorphism invariant algebra.

While Colombeau in his ``Elementary Introduction'' \cite{c2} chose
to embed distributions via convolution with a mollifier, i.e., 
(cf.\ also (\ref{iota}) above)
\beq\label{iotaC}
        \iota^C(w)(\vphi,x)\,:=\,\langle w,\vphi(.-x)\rangle\,,
\eeq
Jel\'\i nek (following in fact earlier ideas of Colombeau presented 
in \cite{c1}) decided to embed distributions by letting them act on 
the test function, i.e.,
\beq\label{iotaJ}
        \iota^J(w)(\vphi,x)\,:=\,\langle w,\vphi\rangle\,.
\eeq
Since both embeddings are simply related by a translation, i.e., 
$\iota^C=T^*\iota^J$ with
\beas
        T:\D(\R^s)\times\R^s&\to&\D(\R^s)\times\R^s\\
        (\vphi,x)&\mapsto&(T_x\vphi,x):=(\vphi(.-x),x)\,,
\eeas
they give rise to equivalent descriptions of virtually every  Colombeau algebra, which 
we call {\em C- resp.\ J-formalism}. In \cite{found}, Section 5 a translation formalism 
allowing to change from one setting to the other at any place of the construction 
was established and used in turn to clarify subtle questions of 
infinite-dimensional calculus. Before giving the actual definitions 
of $\G^d$ we briefly comment on these issues.
Jel\'\i nek uses \cite{Y} as main reference while the presentation
of \cite{found} and \cite{vim} is based upon the more convenient calculus 
of \cite{km}. The basic idea of the latter is that a map $f:E\to F$ between locally
convex spaces is smooth if it transports smooth curves in $E$ to smooth curves in $F$,
where the notion of smooth curves is straightforward (via limits of difference 
quotients).
This notion of smoothness in general is weaker than Silva-differentiability but coincides
with the latter on all the spaces used in the construction of Colombeau algebras.
Moreover, it displays the following decisive advantage in applications to partial 
differential equations: if one is to construct a generalized solution to a 
nonlinear singular equation this is done componentwise, i.e., for fixed $\vphi$. 
Smoothness of the respective solution in $\vphi$ is then guaranteed already by classical 
theorems on smooth dependence of solutions on parameters.

\section{The algebra $\G^d(\Om)$}  

We now give a brief description of Jel\'\i nek's algebra $\G^d$: contrary 
to the original presentation using the $C$-formalism for its better 
familiarity (however, omitting the superscript $C$ from now on). For a comparison 
of the respective features of the two formalisms we refer to the table in 
\cite{found}, Section 5. Apart from closing
a gap in the construction of \cite{Jel} the presentation in \cite{found} supplies
those parts of the resp.\ arguments which have not been included in \cite{Jel}.
This applies, in particular, to the questions of smoothness and stability of 
${\cal E}_M$, $\NN$ w.r.t.\ differentiation and the fact that transformed test objects
are not defined on the whole of $(0,1]\times\Om$ in general.

Forced by the choice of the embedding (\ref{iotaC}) we define the basic space to 
be 
\[
        \E^d(\Om)\,:=\,\CC^\infty(U(\Om))\,.
\]
Partial derivatives on $\E^d(\Om)$---which will become
the derivatives in the algebra---in the C-formalism are simply defined by
\beq\label{DiC}
        D^d_i:\ \E^d(\Om)\to\E^d(\Om)\qquad D^d_i=\pa_i\,.
\eeq

Recall that test objects have to depend on $\eps$ and $x$, in particular are chosen
to be smooth, bounded paths $\phi:(0,1]\times\Om\to\A_0(\R^s)$ (resp.\ $\A_q(\R^s)$). 
Denoting their
space by $\CC^\infty_b((0,1]\times\Om,\A_0(\R^s))$ we are able to formulate the
tests for moderateness and negligibility.

\bdef\label{jtest}
\begin{itemize}
\item[(i)]
$R\in\ce^d(\Om)$ is called {\rm moderate} if\\
$        \forall K\subset\subset\Om\ \forall\al\in\N_0^s\ \exists N\in\N\ 
        \forall \phi\in\cc^\infty_b((0,1]\times\Om,\ca_0(\R^s)):$
\[
        \sup_{x\in K}|\pa^\al(R(S_\eps\phi(\eps,x),x))|
        =O(\eps^{-N})\mbox{ as }\eps\to0\,.
\]
The set of all moderate elements $R\in\ce^d(\Om)$ will be denoted by $\ce^d_M(\Om)$.
\item[(ii)]
$R\in\ce^d_M(\Om)$ is called {\rm negligible} if\\
$\forall K\subset\subset\Om\ \forall\al\in\N_0^s\ \forall n\in\N\ 
        \exists q\in\N\ 
        \forall \phi\in\cc^\infty_b((0,1]\times\Om,\ca_q(\R^s)):$
\[
        \sup_{x\in K}|\pa^\al(R(S_\eps\phi(\eps,x),x))|
        =O(\eps^s)\mbox{ as }\eps\to0\,.
\]
The set of all negligible elements $R\in\ce^d(\Om)$
will be denoted by $\cn^d(\Om)$.
\end{itemize}
\ethi

The key ingredients in proving diffeomorphism invariance as well as stability
with respect to derivatives (i.e., that the $x$-derivative of a moderate resp.\
negligible function again is moderate resp.\ negligible; this becomes a
peculiar issue due to the additional $x$-dependence of $\phi$) are several 
equivalent formulations of the tests given above.

To settle the question of stability w.r.t. differentiation Jel\'\i nek
introduced an alternate, yet equivalent form of tests involving differentials 
of $R$ with respect to the test function-slot denoted by $\rmd_1$. (\cite{Jel}, 
Th.\ 17, resp.\ Th.\ 18, $(2^\circ)\Leftrightarrow(3^\circ)$). 
We only formulate the respective test for
moderateness (the case of negligibility being analogous) and refer to the original
for the ingenious proofs. We presume that the author was completely aware
of the role Ths.\ 17 and 18 had to play in this respect 
yet for some reasons he decided not to address this issue.

\bthm
$R\in\ce^d(\Om)$ is a member of
$\ce^d_M(\Om)$ if and only if the following condition is satisfied:
\beas\forall K\subset\subset\Om\ \forall\al\in\N_0^s\ 
   \forall k\in\N_0\ \exists N\in\N
   \ \forall B\,(\mbox{bounded})\,\subseteq \cd(\R^s):\eeas
$$\pa^\al\rmd_1^k(R\circ S^{(\eps)})(\vphi,x)(\psi_1,\dots,\psi_k)=O(\eps^{-N})
   \qquad\qquad (\eps\to0)$$
uniformly for $x\in K$, $\vphi\in B\cap\ca_0(\R^s)$,
   $\psi_1,\dots,\psi_k\in B\cap\ca_{00}(\R^s)$.
\ethi
Here $\ca_{00}(\R^s)=\{\vphi\in\D(\R^s): \int\vphi=0\}$ 
denotes the tangent space of $\A_0$ and the operator
$S^{(\eps)}$ is derived from $S_\eps$ by
\beas
        S^{(\eps)}:\cd(\R^s)\times\R^s&\to&\cd(\R^s)\times\R^s\\
        (\vphi,x)&\mapsto&(S_\eps\vphi,x)=(\frac{1}{\eps^s}\,\vphi(\frac{.}{\eps}),x)\,.
\eeas
Finally to prove that $\pa_iR$ is moderate if $R$ was, observe that 
$\pa_i(R\circ S^{(\eps)})=(\pa_iR)\circ S^{(\eps)}$. Then the claim follows from
$$
\pa^\al\rmd_1^k((\pa_iR)\circ S^{(\eps)})
  =\pa^\al\rmd_1^k\pa_i(R\circ S^{(\eps)})
  =\pa^{\al+e_i}\rmd_1^k(R\circ S^{(\eps)})\,.
$$

We now turn to the central issue of diffeomorphism invariance. 
First we present a heuristical calculation which clearly shows which path has
to be pursued. Suppose we want to prove moderateness of $\hat\mu R$. Given
$\tilde\phi\in\CC^\infty((0,1]\times\tilde\Om,\A_0(\R^s))$ then we would 
have to estimate
\beas
   \hat\mu(R)(S_\eps\tilde\phi(\eps,\tilde x),\tilde x)
   &=&R(\bar\mu(S_\eps\tilde\phi(\eps,\tilde x),\tilde x))\,=\,
   R(\bar\mu S^{(\eps)}(\tilde\phi(\eps,\tilde x),\tilde x))\\
   &=&R(S^{(\eps)}(S^{(\eps)})^{-1}\bar\mu S^{(\eps)}(\tilde\phi(\eps,\tilde x),
   \tilde x))\,.
\eeas
Hence we would need $R$ to pass a test for moderateness w.r.t. test objects
of the form (denoting by ${\mathrm pr}_1$ the projection to the fist component)
\beas
   \phi(\eps,x)&=&
   {\mathrm pr}_1(S^{(\eps)})^{-1}\bar\mu S^{(\eps)}(\tilde\phi(\eps,\tilde x))\\
   &=&\ti\phi(\eps,\mu^{-1}x)\left(\frac{\mu^{-1}(\eps\xi+x)-\mu^{-1}x}
   {\eps}\right)\cdot|\det D\mu^{-1}(\eps\xi+x)|\,.
\eeas
But unfortunately $\phi(\eps,x)\not\in\CC^\infty((0,1]\times\Om,\A_0(\R^s))$
since it is only defined if $\xi\in\frac{\Om-x}{\eps}$, whereas we want
$\xi\mapsto\phi(\eps,x)(\xi)$ to be a test function on the whole of $\R^s$. 

However, $\phi(\eps,x)$ belongs to a class of test objects providing an apparently 
weaker, yet, as it finally turns out, equivalent test. More
precisely, from \cite{found}, Th.\ 10.5 we have that $R\in\E^d(\Om)$ is moderate
if and only if it fulfills the following condition
(Z) \footnote{Condition (Z) is one of 6 tests proved to be equivalent
in \cite{found}, Theorem 10.5 called Theorem A--Z there. Note, however,
that this neither indicates that the authors of \cite{found}
originally intended to give 26 equivalent tests, nor that (Z) for some
mysterious reason was considered to be the ultimate condition;
rather it was invented during a train ride returning from a
workshop in Novi Sad to Vienna and ``Z'' simply 
stands for ``Zug'' which is the German word for train.}
\[
        \forall K\subset\subset\Om\ \forall\al\in\N_0^s\ \exists N\in\N
        \ \forall\phi
        :D\to\ca_0(\R^s))\ \mbox{($D,\phi$ as described below})
\]
\[
        \exists C>0\ \exists \eta>0\ \forall\eps\,(0<\eps<\eta)
        \ \forall x\in K:\ (\eps,x)\in D \mbox{ and }
\]\[
        |\pa^\al(R(S_\eps\phi(\eps,x),x))|\leq C\eps^{-N}\,,
\]
where $D\subseteq(0,1]\times\Om$ and for
$D,\phi$ the following holds: For each $L\subset\subset\Om$ there
exists $\eps_0$ and a subset $U$
of $D$ which is open in $(0,1]\times\Om$ such that
\begin{itemize}
\item[(1)] 
          $(0,\eps_0]\times L\subseteq U(\subseteq D)$
          and $\phi$ is smooth on $U$, and
\item[(2)] for all $\bet\in\N_0^s$,
          $\{\pa^\bet\phi(\eps,x)\mid0<\eps\le\eps_0,\ x\in L\}$ is bounded in $
        \cd(\R^s)$.
\end{itemize}
Now diffeomorphism invariance of the notion of moderateness is 
established by the following
\bthm\label{T6} Let $\mu:\ti\Om\to\Om$ be a diffeomorphism and 
$\ti\phi\in\cc^\infty_b((0,1]\times\ti\Om,\ca_0(\R^s))$.
Define $D(\subseteq(0,1]\times\Om$) by
\[D:=\{(\eps,x)\in(0,1]\times\Om\mid(\ti\phi(\eps,\mu^{-1}x),\mu^{-1}x)
\in U_\eps(\ti\Om)\}\,.
\]
For $(\eps,x)\in D$, set
$$\phi(\eps,x)(\xi):=
   \ti\phi(\eps,\mu^{-1}x)\left(\frac{\mu^{-1}(\eps\xi+x)-\mu^{-1}x}
   {\eps}\right)\cdot|\det D\mu^{-1}(\eps\xi+x)|.$$
Then $\phi$ satisfies the requirements specified for test objects
in condition (Z).
\ethi
In some more detail assume $R$ to be moderate. We show that $\hat\mu R$ passes the test
used in Definition \ref{jtest} (i). Indeed given $\tilde K\subset\subset\tilde\Om$,
$\al\in\N^s_0$ and $\tilde\phi\in\CC^\infty_b((0,1]\times\tilde\Om,\A_0(\R^s))$
define $\phi$ as in the preceding theorem. Then by the chain rule
\beas
     \pa_{\ti x}^\al
     \left((\hat\mu R)(S_\eps\ti\phi(\eps,\ti x),\ti x)\right)&=&
      \pa_{\ti x}^\al
      (R(S_\eps\phi(\eps,\mu\ti x),\mu\ti x))\\
      &=&\sum\limits_{|\,\bet|\le|\al|}
      \pa_{x}^\bet
      (R(S_\eps\phi(\eps,x),x))\Big|_{x=\mu\ti x}
      \cdot g_\bet(\ti x)\,,
\eeas
where each $g_\beta$ is bounded on $\tilde K$. Since $R$ satisfies condition (Z) 
the claim follows.

However, matters become more complicated in the case of negligibility. 
First note that the resp.\ test objects take values in $\A_q$
($q>0$) which is not a diffeomorphism invariant property. The way out
is provided by the re-introduction of asymptotically vanishing
moments (cf.\ \ref{ColM}) into the theory by building up (another) equivalent test
using this notion. Indeed Jel\'\i nek (\cite{Jel}, 18 (4$^\circ$)) has formulated 
such a condition which, however, unfortunately is {\em not equivalent} to
the notion of negligibility as defined above. Moreover, this condition 
is so strong that even $\iota(x^2)\not=(\iota(x))^2$, hence the property 
of $\iota$ being an algebra homomorphism
on $\CC^\infty$ is lost. However, in \cite{found}, Section 7
this flaw was removed, namely by demanding also {\em all derivatives} 
of the test objects to have asymptotically vanishing moments. 
More precisely we say that a test object 
$\phi\in\CC^\infty_b((0,1]\times\Om,\A_0(\Om))$ is of
type $[A_l^\infty]_{K,q}$ if on (a given) $K\subset\subset\Om$ 
\[
        \forall\beta\in\N^s_0\ 1\leq|\beta|\leq q\ \forall\ga\in\N^s_0
        \ \sup_{x\in K}|\int \xi^\be\pa^\ga\phi(\eps,x)(\xi)\,dx|=O(\eps^q)\,.
\]
Then we have

\bthm
$R\in\E^d(\Om)$ is negligible if and only if $\forall K\subset\subset\Om\ 
\forall\al\in\N_0^s$ $\forall n\in\N$ $\exists q$ such that $\forall
\phi$ of type $[A_l^\infty]_{K,q}$:
\[ 
        \sup_{x\in K}|\pa^\al R(S_\eps\phi(\eps,x),x)|=O(\eps^n)\,.
\]
\ethi
Diffeomorphism invariance of the notion of negligibility is now established 
by the above theorem in conjunction
with an analog of Theorem \ref{T6}, as well as an analog of condition (Z)
above. Finally we define our main object of desire.
\bdef
The {\em diffeomorphism invariant Colombeau algebra on $\Om$} is defined as the quotient
\[
        \G^d(\Om):=\E^d_M(\Om)\,/\,\NN^d(\Om)\,.
\]
\ethi

Summing up we have constructed a differential algebra $\G^d(\Om)$ in a
diffeomorphism invariant way, in particular allowing for a 
diffeomorphism invariant embedding of distributions. Moreover, $\G^d(\Om)$
(as usual) is a fine sheaf of differential algebras (\cite{found}, Section 8).

We finally turn to the issue of commutativity of the embedding with partial
derivatives in the algebra. This will guarantee that $\G^d(\Om)$
indeed possesses all the favorable properties of a Colombeau algebra.
To this end it is useful to change to the J-formalism. Recall from (\ref{DiC})
that derivatives
in the C-formalism are just given by partial derivatives. Using the translation
formalism of \cite{found} we derive that in the J-formalism, i.e., on
$\E^J(\Om):=(T^{-1})^*(\E^d(\Om))$ we have
\beq\label{DiJ}
        D_i^J=(T^{-1})^*\circ\pa_i\circ T^*
        \mbox{ i.e., }
        (D_i^JR)(\vphi,x)=-(d_1R(\vphi,x))(\pa_i\vphi)+(\pa_iR)(\vphi,x)\,.
\eeq
We now see immediately that if $F\in\D'(\Om)$ then $\iota^J(F)$ (cf.\ (\ref{iotaJ}))
is independent of $x$ hence the second term in (\ref{DiJ}) vanishes. Moreover, since
$\iota^J(F)$ is linear in $\vphi$, $-d_1(\iota^J(F))(\vphi,x)(\pa_i\vphi)=
\langle F,-\pa_i\vphi\ra$ which is exactly the $\iota^J$-image of $\pa_iF$.


\begin{thebibliography}{99}
\bibitem[{Col}82]{c0}
{Colombeau, J.~F.}
\newblock {\em Differential Calculus and Holomorphy, Real and Complex Analysis
  in Locally Convex Spaces}.
\newblock North Holland, Amsterdam, 1982.

\bibitem[{Col}84]{c1}
{Colombeau, J.~F.}
\newblock {\em New Generalized Functions and Multiplication of Distributions}.
\newblock North Holland, Amsterdam, 1984.

\bibitem[{Col}85]{c2}
{Colombeau, J.~F.}
\newblock {\em Elementary Introduction to New Generalized Functions}.
\newblock North Holland, Amsterdam, 1985.

\bibitem[{Col}94]{CM}
{Colombeau, J.~F., Meril, A.}
\newblock Generalized functions and multiplication of distributions on
  {${\mathcal C}^\infty$} manifolds.
\newblock {\em J.~Math.~Anal.~Appl.}, {\bf 186}:357--364, 1994.

\bibitem[{Gro}99]{vim}
{Grosser, M., Kunzinger, M., Steinbauer, R., Vickers, J.}
\newblock A global theory of algebras of generalized functions.
\newblock {\em Preprint (available electronically at {\tt
  http://arXiv.org/abs/math.FA/9912216})}, 1999.

\bibitem[{Gro}01]{found}
{Grosser, M., Farkas, E., Kunzinger, M., Steinbauer, R.}
On the foundations of nonlinear generalized functions {I}, {II}.
{\em Mem.~Am.~Math.~Soc., to appear (available electronically at {\tt
  http:}\,{\tt //arXiv.org/abs/math.FA/9912214, 9912215})}, 2001.

\bibitem[{Jel}99]{Jel}
{Jel\'\i nek, J.}
\newblock An intrinsic definition of the {C}olombeau generalized functions.
\newblock {\em Comment.~Math.~Univ.~Carolinae}, {\bf 40}:71--95, 1999.

\bibitem[{Kri}97]{km}
{Kriegl, A., Michor, P.~W.}
\newblock {\em The Convenient Setting of Global Analysis}, volume~{\bf 53} of
  {\em Math.~Surveys Monogr.}
\newblock Amer.~Math.~Soc., Providence, RI, 1997.

\bibitem[{Kun}96]{kunzDISS}
{Kunzinger, M.}
\newblock {\em Lie Transformation Groups in {C}olombeau Algebras}.
\newblock PhD thesis, University of Vienna, 1996.

\bibitem[{Kun}99]{geo2}
{Kunzinger, M., Steinbauer, R.}
\newblock A rigorous solution concept for geodesic and geodesic deviation
  equations in impulsive gravitational waves.
\newblock {\em J.~Math.~Phys.}, {\bf 40}:1479--1489, 1999.

\bibitem[{Kun}01]{ndg}
{Kunzinger, M., Steinbauer, R.}
Nonlinear distributional geometry.
{\em Preprint (available electronically at {\tt http://}\,
{\tt arXiv.org/abs/math.FA/0102019})}, 2001.

\bibitem[{Ned}98]{NPS}
{Nedeljkov, M., Pilipovi\'c, S., Scarpal\'ezos, D.}
\newblock {\em The linear theory of Colombeau generalized functions},
  volume~{\bf 385} of {\em Pitman Research Notes in Mathematics}.
\newblock Longman, Harlow, 1998.

\bibitem[{Vic}99]{invc}
{Vickers, J.~A., Wilson, J.~P.}
\newblock Invariance of the distributional curvature of the cone under smooth
  diffeomorphisms.
\newblock {\em Class.~Quantum.~Grav.}, {\bf 16}:579--588, 1999.
\bibitem[{Yam}74]{Y}
{Yamamuro, S.}
\newblock {\em Differential Calculus in Topological Linear Spaces}, volume~{\bf
  374} of {\em Lecture Notes in Mathematics}.
\newblock Springer, New York, 1974.   
\end{thebibliography}
\end{document}